%% file: Main.tex
\documentclass{article}

\usepackage{arxiv}

\usepackage{graphicx}
\usepackage[utf8]{inputenc} 
\usepackage[T1]{fontenc}    
\usepackage[colorlinks,allcolors=blue]{hyperref}       
\usepackage{url}            
\usepackage{booktabs}       
\usepackage{amsfonts}       
\usepackage{nicefrac}       
\usepackage{microtype}      
\usepackage{lipsum}		
\usepackage{graphicx}
\usepackage{natbib}
\usepackage{doi}
\usepackage{tabularx,booktabs}
\usepackage{algorithm}
\usepackage{algpseudocode} 
\usepackage{multirow}
\usepackage{xcolor}
\usepackage{multicol}
\usepackage{subcaption}
\usepackage{tikz} 
\usepackage{amsmath}
\usepackage{textcomp}
\usepackage{libertine}
\usepackage{comment}
\usepackage{amssymb}

\algnewcommand\algorithmicforeach{\textbf{for each}}
\algnewcommand\algorithmicdoparallel{\textbf{do in parallel}}
\algdef{S}[FOR]{ForEach}[1]{\algorithmicforeach\ #1\ \algorithmicdo}
\algdef{S}[FOR]{ForEachParallel}[1]{\algorithmicforeach\ #1\ \algorithmicdoparallel}
\algnewcommand{\sIf}[2]{
  \State \algorithmicif\ #1\ \algorithmicthen\ #2}

\graphicspath{graphics}
\usepackage{xspace}

\title{Penalty Weights in QUBO formulations: Permutation Problems\thanks{Please cite: \protect\url{https://doi.org/10.1007/978-3-031-04148-8_11}}}

\author{ \href{https://orcid.org/0000-0003-0854-4777}{\includegraphics[scale=0.06]{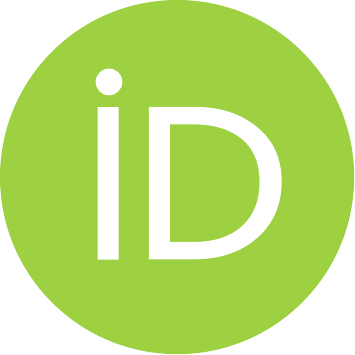}\hspace{1mm}Mayowa Ayodele}\\
	Fujitsu Research of Europe Ltd.\\
	The Urban Building, 3-9 Albert Street\\
	Slough, United Kingdom,  SL1 2BE\\
	\texttt{mayowa.ayodele@fujitsu.com} \\
}

\renewcommand{\shorttitle}{Penalty Weights in QUBO formulations: Permutation Problems}

\begin{document}
\maketitle

\begin{abstract}
Optimisation algorithms designed to work on quantum computers or other specialised hardware have been of research interest in recent years. Commercial solvers that use quantum or quantum-inspired methods, such as Fujitsu’s Digital Annealer (DA) and D-wave’s Quantum Annealer, can solve optimisation problems faster than algorithms implemented on general purpose computers. However, they can only optimise problems that are in binary and quadratic form. Quadratic Unconstrained Binary Optimisation (QUBO) is therefore a common formulation used by these solvers.

There are many combinatorial optimisation problems that are naturally represented as permutations e.g. travelling salesman problem. Encoding permutation problems using binary variables however presents some challenges. Many QUBO solvers are single flip solvers, it is therefore possible to generate solutions that cannot be decoded to a valid permutation. To create bias towards generating feasible solutions, we use penalty weights. The process of setting static penalty weights for various types of problems is not trivial. This is because values that are too small will lead to infeasible solutions being returned by the solver while values that are too large may lead to slower convergence. In this study, we explore some methods of setting penalty weights within the context of QUBO formulations. We propose new static methods of calculating penalty weights which lead to more promising results than existing methods.

\end{abstract}


\input{sections_final}


\bibliographystyle{unsrtnat}
\bibliography{References}


\end{document}

%% file: sections_final.tex
\section{Background}
Permutation problems are well studied combinatorial optimisation problems in nature inspired computing. They have many real-world applications especially in planning and logistics. Some of the most frequently studied permutation problems in literature are the well-known Travelling Salesman Problem (TSP) and Quadratic Assignment Problem (QAP). Since these problems are NP-hard, heuristics and meta-heuristics have been proposed for solving them. 

Several classes of algorithms have been applied to problems naturally represented as permutations e.g. Estimation of Distribution Algorithm \cite{arza2020kernels}, Iterated Local Search and Differential Evolution Algorithm \cite{santucci2015algebraic}. Quantum-inspired algorithms such as the Digital Annealer (DA) has also been shown to present more promising performance on the QAP when compared to CPLEX and QBSolve \cite{matsubara2020digital}. The DA was particularly shown to be up to three or four orders of magnitude faster than CPLEX on QAP instances and maximum cut instances.

Quantum and quantum-inspired methods have been of research interest in recent years. This is because they are able to exploit the use of specialised hardware to solve optimisation problems much quicker than classical algorithms implemented on general purpose machines \cite{aramon2019physics}. It is common to formulate combinatorial optimisation problems such as permutation problems in quadratic and binary form such that algorithms that use specialised hardware including (but not limited to) quantum devices can be used to solve them. In recent years, Quadratic Unconstrained Binary Optimisation (QUBO) has become a unifying model for representing many combinatorial optimisation problems \cite{verma2020penalty}. QUBO (or the equivalent Ising) formulations of common combinatorial optimisation problems are presented in \cite{lucas2014ising}. As the name depicts, QUBO problems are unconstrained, quadratic and of binary form. Since the representation only supports binary information, the natural representation of permutation problems can therefore not be used. While some classical optimisation algorithms such as genetic algorithms use permutation representation \cite{hussain2017genetic} to solve QAP and/or TSP, other algorithms require alternative representations e.g. random keys \cite{ayodele2016rk}, factoradics \cite{regnier2014factoradic}, binary \cite{baluja1994population}, matrix representation \cite{larranaga1999genetic} or two-way one-hot \cite{lucas2014ising}. The two-way one-hot representation \cite{glover2018tutorial,lucas2014ising,matsubara2020digital}, also known as permutation matrix \cite{Birdal_2021_CVPR}, is often used in QUBO formulations of permutation problems and is used in this study. In this representation, a substring of bits is used to represent each entity (e.g a location in TSP). In each substring, only one bit can be turned on and in the entire solution, the bit turned on in each substring must be unique. This representation ensures that the mutual exclusivity constraint of permutation problems is respected. QUBO solvers such as Path Relinking method used in \cite{verma2020penalty} and the first generation DA in \cite{aramon2019physics} are single flip solvers and are therefore not able to preserve two-way one-hot validity. To ensure that the problem to be solved is in `unconstrained' form, penalty weights are applied to combine the cost function (unconstrained objective function) with the constraint function. Solutions are penalised by the magnitude of violation of the constraint(s). 

Setting penalty weights is not a trivial task. Values that are too large make the search too difficult for the solver as the penalty terms overwhelm the original objective function \cite{verma2020penalty}. Penalty weights that are too small are highly undesirable as infeasible solutions will displace feasible solutions in the search, causing the solver to return infeasible solution(s). Penalty weights are however not unique and there are often a range of values that work well \cite{glover2019quantum}. The primary objective of setting the penalty weights is often to ensure that the optimal solution of the QUBO is the optimal solution of the original constrained problem. However, it is also important that these values are not too large. 
 
In literature, there are many approaches of setting penalty weights for QUBOs. The value can be set by domain experts \cite{glover2019quantum}. A common approach is to derive the penalty weights empirically using methods that increase the weights until feasibility is achieved \cite{rosenberg2016solving}. This approach is however computationally intensive as a full run of the algorithm and analysis of results is required each time until the right value is reached.

Another common approach is to set the penalty weight to a value greater than the largest possible objective value. Deriving the range (upper bound and lower bound) of the objective function is often problem specific. These values can also be too large. Although there are some general methods of determining the bounds of a QUBO, it is however often the case that methods with less computational complexity lead to values that are too large while methods that can provide better bounds are often computationally expensive \cite{boros2008max}. Moreover, while better bounds can lead to smaller but valid penalty weights (i.e guarantees feasibility of the optimal), the penalty values are still often larger than desired.  In \cite{cseker2020digital}, the performance of the DA was analysed using different penalty weights that are fractions of the range of the objective function. The range was derived using problem specific information, the authors found much smaller values can lead to better performance of the DA.

Furthermore, the maximum coefficient of the QUBO has been used as penalty weights when solving problems like the TSP. The idea behind this is that, if a TSP solution is penalised by the maximum distance between any two cities, feasibility of the optimal solution can always be achieved. In \cite{takehara2019multiple}, the authors used a multi-trial approach. The values used were within a range defined using the minimum and the maximum distance between any two locations. Similarly, in \cite{goh2020hybrid}, fractions of the maximum distance were used to derive penalty weights for the TSP. This is an example of a scenario where values much lower than the full range of the objective function can be valid. This approach may however not generalise to other problems \cite{verma2020penalty}.

In \cite{verma2020penalty}, a pre-processing method that can be used to generate penalties for equality constraints within the context of single flip QUBO solvers was presented. The authors measure the maximum change in objective function that can be obtained as a result of any single flip in a solution. This method, which is referred to as VLM in this study, is explained in more details in Section \ref{sec:pre}.

Examples of QUBO (or Ising) solvers that use specialised hardware are D-wave's Quantum Annealer \cite{johnson2011quantum} and Fujitsu's DA \cite{aramon2019physics}. The Quantum Annealer uses Quantum Processing Units (QPU) to achieve its speed up while the DA uses application-specific CMOS hardware. In this study, we analyse the effect of different methods of generating penalty weights for permutation problems (TSP and QAP). Initial experiments are based on a CPU implementation of the DA Algorithm presented in \cite{aramon2019physics}. Further analysis of the effect of penalty weights are done using the third generation DA (DA3) as the QUBO solver \cite{da3}. 

The rest of this paper is structured as follows. Section \ref{sec:pe} presents the permutation formulations and QUBO formulations of the TSP and QAP used in this study. Section \ref{sec:pre} presents the methods of generating penalty weights. Section \ref{sec:da} presents a description of the DA Algorithm. Experimental Settings are presented in Section \ref{sec:ex}. Analyses of results and conclusion are presented in Sections \ref{sec:res} and \ref{sec:con} respectively.

\section{Permutation Problems}
\label{sec:pe}
A valid permutations is described as $\sigma = \left \{ \sigma_1,\ldots, \sigma_n \right \}$, where  $\sigma_i \neq \sigma_j\ \forall\ i \neq j$. In general, QUBO problems can be defined as follows:

\begin{equation}
\label{eq:qubo}
  E(x) = x^TQx + k \;,
\end{equation}

where $Q$ and $k$ are $m \times m$ QUBO matrix and constant term, the solution $x=(x_1, \dots,x_m)$ is an $m$-dimensional vector, and $E(x)$ is the energy (or fitness) of $x$. To formulate permutation problems as QUBO, it is important for the problem to be formulated as zeros and ones as these are the only values supported by QUBO solvers. Some of the well-known approaches of transforming integer values to zeros and ones are binary, gray code and one-hot. Within the context of QUBOs, the two-way one-hot (permutation matrix) encoding is often used to represent permutations and will be used in this study.

The energy (fitness function) of the problem is shown in Eq.~(\ref{eq:perm}) where $c(x)$ and $g(x)$ are respectively cost and constraint functions while $\alpha$ is the penalty weight. The generic cost and constraint functions of permutation problems are presented in Eqs. (\ref{eq:quboperm1}) and (\ref{eq:penfun}) \cite{goh2020hybrid}.

\begin{equation}
\label{eq:perm}
  \text{minimise } E(x) = c\left ( x \right ) + \alpha \times g\left ( x \right )
\end{equation}

\begin{equation}
\label{eq:quboperm1}
     c(x) = \sum_{i=1}^{n}\sum_{j=1}^{n}\sum_{k=1}^{n}\sum_{l=1}^{n}x_{i,k}Q_{i,k,j,l}x_{j,l} 
\end{equation} 

\begin{equation}
\label{eq:penfun}
g(x) = \sum_{i=1}^{n}\left ( 1- \sum_{k=1}^{n} x_{i,k}\right )^2 + \sum_{k=1}^{n}\left ( 1- \sum_{i=1}^{n} x_{i,k}\right )^2 \;.
\end{equation}

Note that the constraint function is designed to ensure that the solutions can be converted to valid permutations i.e. penalise solutions that do not satisfy,

\begin{align}
 \label{eq:quboperm2}
    \sum_{k=1}^{n}x_{i,k} = 1\ \forall\ i \in \left \{  1,\ldots,n\right \}, \quad & 
    \sum_{i=1}^{n}x_{i,k} = 1\ \forall\ k \in \left \{  1,\ldots,n\right \} \;.
\end{align}

In Eqs. \eqref{eq:quboperm1} - \eqref{eq:quboperm2}, binary variable $x_{i,k}$ indicates whether an object $i$ is assigned to position $k$ or not. We however note that $x$ is solved as a vector of size $m=n^2$ rather than a $n \times n$ matrix. Also, while $Q_{i,k,j,l}$ in Eq.~\eqref{eq:quboperm1} is the QUBO coefficient that captures the relationship between an object $i$ being in position $k$ and an object $j$ being in position $l$. In the rest of this paper, QUBO matrices $C$ and $G$ representing the cost or constraint functions are presented as $m \times m$ matrices.

\subsection{Quadratic Assignment Problem}
The QAP can be described as the problem of assigning a set of $n$ facilities to a set of $n$ locations. For each pair of locations, a distance is specified. For each pair of facilities, a flow (or weight) is specified. The aim is to assign each facility to a unique location such that the sum of the products between flows and distances is minimised. 

Formally, the QAP consists of two $n \times n$ input matrices $H=[h_{i,j}]$ and $D=[d_{k,l}]$, where $h_{i,j}$ is the flow between facilities $i$ and $j$, and $d_{k,l}$ is the distance between locations $k$ and $l$, the solution to the QAP is a permutation $\sigma = (\sigma_1, \ldots ,\sigma_n )$ where $\sigma_i$ represents the location that facility $i$ is assigned to. The objective function (total cost) is formally defined as follows

\begin{equation}
    \text{minimise}\ f(\sigma) = \sum_{i=1}^{n}\sum_{j=1}^{n}h_{i,j} \times d_{\sigma_i, \sigma_j}\;.
\end{equation}

We aim to solve Eq.~(\ref{eq:perm}) where the cost function $c(x)$ of the QUBO representing the QAP is presented in Eq.~(\ref{eq:cqap}) and the constraint function $g(x)$ is the same as Eq.~(\ref{eq:penfun}).

\begin{equation}
\label{eq:cqap}
  c\left ( x \right ) = \sum_{i=1}^{n}\sum_{j=1}^{n}\sum_{k=1}^{n}\sum_{l=1}^{n}h_{i,j}d_{k,l}x_{i,k}x_{j,l}\;.
\end{equation}

Any solution $x$ can be encoded with $n^2$ variables, when $x$ is presented in vector format.

\subsection{Travelling Salesman Problem}
    The TSP consists of $n$ locations and a matrix $d$ representing distances between any two locations. The aim of the TSP is to minimise the distance travelled while visiting each location exactly once and returning to the location of origin. Given that $\sigma_i$ is used to denote the $i^{th}$ city and $d_{\sigma_{i-1},\sigma_i}$ is the distance between $\sigma_i$ and $\sigma_{i-1}$. The solution to the TSP is a permutation $\sigma = \left \{ \sigma_1,\ldots ,\sigma_n \right \}$ where each $\sigma_i\ (i= 1,\ldots ,n)$ represents the $i^{th}$ location to visit. The TSP is formally defined as

\begin{equation}
\label{eq:tsp}
    \text{minimise } f\left ( \sigma \right ) = \sum_{i=2}^{n} d_{\sigma_{i-1}, \sigma_i} + d_{\sigma_n, \sigma_1}\;.
\end{equation}

We aim to solve Eq.~(\ref{eq:perm}) where the cost function, $c(x)$ of the QUBO representing the TSP is presented in Eq.~(\ref{eq:qtsp}) and the constraint function $g(x)$ is the same as Eq.~(\ref{eq:penfun}).

\begin{equation}
\label{eq:qtsp}
c\left ( x \right ) = \sum_{(l,i) \in E}d_{l,i} \sum_{k=1}^{n}x_{l,k}x_{i,k+1} \;.
\end{equation}

The TSP instances used in this work are symmetric, we can therefore fix the first city, reducing the size of $x$ to $(n-1)^2$ \cite{lucas2014ising}. 

Note that QUBOs for these problems can be generated using packages such as PyQUBO \cite{zaman2021pyqubo}. QUBOs are generated in this study using method in \cite{moraglio_georgescu}.

\section{Penalty Weights}
\label{sec:pre}
The aim of this study is to derive methods of setting $\alpha$ in Eq.~(\ref{eq:perm}) such that the optimal solution to the penalised objective function is the optimal solution of the original constrained problem. We do this without problem specific knowledge but use information captured in the QUBO matrices representing the cost and constraint functions. This is shown in Eq.~(\ref{eq:ca}), $c\left ( y \right )$ is used to denote the cost function of the optimal solution $y$. $S$ is the solution space of infeasible solutions. Note that $g(x)$ produces a non-negative value, $g(x)=0$ if the solutions are feasible but $g(x)>0$ for infeasible solutions. The value of $g(x)$ increases according to the degree of constraint violation.

\begin{equation}
\label{eq:ca}
c\left ( y \right ) <  c\left(x\right) + \alpha \times g\left(x\right)\ \forall\ x \in S\;.
\end{equation}

Eq.~(\ref{eq:ca}) implies that a valid penalty weight $\alpha$ is one that satisfies Eq.~(\ref{eq:wv})

\begin{equation}
\label{eq:wv}
\alpha > \max_{x \in S}  \left ( \frac{ c\left ( y \right ) - c\left(x\right)}{g\left(x\right)} \right )\;.
\end{equation}

In the rest of this study, $C$ and $G$ are used to denote the QUBO matrices representing $c(x)$ and $g(x)$ respectively. Note that $Q$, which is the QUBO matrix optimised by the solver, can be derived by aggregating the matrices (i.e. $Q = C + \alpha \times G$), where $\alpha \geq 1$. The methods of generating penalty weights used in this study are described as follow:

\textbf{UB:} A common method of setting penalty weights is based on the Upper Bound (UB) of the objective function. The UB of $C$ used in this study is presented in Eq.~\eqref{equb}. This is a valid upper bound for problems with all positive QUBO coefficients. We note that a solution consisting of all 1s is an infeasible solution but it gives an estimate of how large the objective function could be. 

\begin{equation}
\label{equb}
\text{UB} = z^TCz, \  z_i = 1\ \forall\ i\ \in [1,n]\;.
\end{equation}

\textbf{MQC:} The Maximum QUBO Coefficient (MQC) which also corresponds to the maximum distance between any two cities in the TSP has been used as penalty weights in previous study \cite{lucas2014ising}. MQC is defined in Eq.~\eqref{eqmqc}.

\begin{equation}
\label{eqmqc}
\text{MQC} = \max_{i=1}^{n}\max_{j=1}^{n}C_{i,j} \;.
\end{equation}

\textbf{VLM:} This is the method proposed by Verma and Lewis \cite{verma2020penalty}. For a 1-flip solver like the DA, any variable $x_i$ can be flipped from 0 to 1 and vice versa at each iteration of the algorithm. VLM focuses on deriving a good estimate for the numerator of Eq.~(\ref{eq:wv}), i.e. $(c(y) - c(x))$. The method estimates the amount of gain/loss in objective function that can be achieved by either turning a bit on or off. They do not consider the denominator (i.e $g(x)$) which the authors recognise will be hard to estimate without complete enumeration. Since $g(x) > 0$ in infeasible solutions, VLM (Eq.~(\ref{eq:wvermaa})) will always be valid. 

\begin{align}
\label{eq:wvermaa}
W^c = \Biggl\{ -C_{i,i} - \sum_{\substack{j=1\\ j \neq i}}^{n}  \min\{C_{i,j}, 0\},\ & C_{i,i}
+ \sum_{\substack{j=1\\ j \neq i}}^{n} \max\{C_{i,j}, 0\}\  \forall\ i \in  \left [ 1,n \right ] \Biggr\} \\
&\alpha = \text{VLM} =\max_{i=1}^{n}W^c_i\;.
\end{align}

\textbf{MOMC:} We propose an amendment to VLM \cite{verma2020penalty}. We refer to the proposed method as the Maximum change in Objective function divided by Minimum Constraint function of infeasible solutions (MOMC). We note that $g(x)$ is not considered in Eq.~(\ref{eq:wvermaa}). VLM can be reduced such that $\alpha$ is still valid, if we know the minimum constraint function ($g(x)$) of any infeasible solution. This can be computed from $G$ by estimating the minimum change in constraint function that is greater than 0 as shown in Eq.~(\ref{eq:gams}).

\begin{align}
\label{eq:gams}
W^g = \Biggl\{ -G_{i,i} - \sum_{\substack{j=1\\ j \neq i}}^{n}  \min\{G_{i,j}, 0\},\ & G_{i,i}
+ \sum_{\substack{j=1\\ j \neq i}}^{n} \max\{G_{i,j}, 0\}\  \forall\ i \in  \left [ 1,n \right ] \Biggr\} \\
\label{eq:gams2}
&\gamma = \min_{\substack{i=1\\W^g_i >0} }^{n}W^g_i \;.
\end{align}

For permutation problems represented as two-way one-hot, $g(x)$ of any solution that is a flip away from a feasible solution is $2$ (i.e. $\gamma = 2$). Method of generation $\alpha$ using the proposed MOMC is presented in Eq.~\ref{eq:momc} where $W^c$ and $\gamma$ are derived as shown in Eqs. \eqref{eq:wvermaa} and  \eqref{eq:gams2}

\begin{align}
\label{eq:momc}
\alpha =\text{MOMC} = \max\left ( 1, \frac{\text{VLM}}{\gamma} \right ) =\max\left ( 1, \frac{\max_{i=1}^{n}W^c_i}{2} \right )\;.
\end{align}

\textbf{MOC:}
We propose another amendment to the VLM method. The method presented here is derived by selecting the Maximum value derived from dividing each change in Objective function with the corresponding change in Constraint function (MOC). MOC captures possible equivalent increase in constraint function as a result of a change in objective function which could be achieved by flipping any bit from 0 to 1 or vice versa. Method of generation $\alpha$ using the proposed MOC is presented in Eq.~\eqref{eq:moc} where $W^c$ and $W^g$ are derived as shown in Eqs. \eqref{eq:wvermaa} and \eqref{eq:gams}

\begin{align}
\label{eq:moc}
\alpha = \text{MOC} = & \max\left ( 1, \max_{\substack{i=1\\W^g_i > 0}}^{n} abs\left (\frac{W^c_i}{W^g_i}\right )\right ) \;.
\end{align}

\section{Digital Annealer}
\label{sec:da}
The DA is a technology designed to solve large scale combinatorial optimisation problems in much shorter time than most classical algorithms.

\begin{algorithm}[h!tb]
\caption{The DA (1st generation) Algorithm}\label{alg:da}
\begin{algorithmic}[1]
\State initial\_state $\leftarrow$ an arbitrary state \label{a1}
\ForEach {run}
\State initialise to initial\_state
\State $E_{\text{offset}}\ \leftarrow\ 0$
\ForEach {iteration }
\State update the temperature if due for temperature update
\ForEach {variable $j$, in parallel }
\State propose a flip using $\Delta E_j - E_{\text{offset}}$
\State if acceptance criteria ($P_j$) is satisfied, record \label{a2}
\EndFor
\If{at least one flip is recorded as meeting the acceptance criteria} 
    \State chose one flip at random from recorded flips
    \State update the state and effective fields, in parallel
    \State $E_{\text{offset}}\ \leftarrow\ 0$
\Else
    \State $E_{\text{offset}} = E_{\text{offset}} +$  offset\_increase\_rate
\EndIf

\EndFor
\EndFor

\end{algorithmic}
\end{algorithm}

The first generation DA is a single flip solver with similar properties as the Simulated Annealing (SA). It is however designed to be more effective than the classical SA algorithm \cite{aramon2019physics}. The standard algorithm of the first generation DA is presented in Alg.~\ref{alg:da}. The DA algorithm exploits the use of specialised hardware such that all neighbouring solutions are explored in parallel and in constant time regardless of the number of neighbours. This approach significantly improves acceptance probabilities of the regular SA algorithm. The DA does not completely evaluate each solution but computes the energy difference resulting from flipping any single bit of the parent solution. Also, the DA uses an escape mechanism to avoid being trapped in local optimal. As shown in the algorithm, $E_{\text{offset}}$ is used to relax acceptance criteria. The $E_{\text{offset}}$ is incremented by a parameter ($\text{offset\_increase\_rate}$) each time no neighbour which satisfies the acceptance criteria is found. The acceptance criteria used in this study is $P_j = \text{exp}(\text{min}(0,-(\Delta E_j - E_{\text{offset}}) /\delta))$ where $P_j$ is the probability of accepting the $j^{th}$ flip. Note that $\Delta E_j$ represents the difference in energy as a result of flipping the $j^{th}$ bit and $\delta$ is the current temperature.

More extensive details of the algorithm can be seen in \cite{aramon2019physics}. The first and second generation DAs were released in May 2018 and December 2018. Both versions were designed to solve optimisation problems that have been formulated as QUBOs. The most recent generation of the DA is the DA3 which is able to find optimal or sub optimal solutions to Binary Quadratic Problems (BQP) of up to 100,000 bits \cite{da3}. BQPs include QUBO but also other binary and quadratic formulations that may have constraints. Note that the algorithm that supports DA3 has been updated to perform better than the current algorithm presented. For simplicity, we use the CPU implementation of Algorithm~\ref{alg:da} in this study to evaluate the quality of solution derived when using different penalty weights. We however also presented some results derived using DA3.

\section{Experimental Setup}
\label{sec:ex}
Problem sets, measure of performance and parameter settings used in this study are presented in this section.

\subsection{Datasets}
In order to compare the behaviour of the DA with different penalty weights, we used common TSP and QAP instances from TSPLIB \cite{reinelt1991tsplib} and QAPLIB \cite{burkard1997qaplib}. The instances used in this study and the corresponding solution sizes $m$ are presented in Table \ref{tb:inst}.

\begin{table}[h]
\centering
\caption{QAP and TSP instances and their solution sizes \label{tb:inst}}
\resizebox{0.8\textwidth}{!}{
\begin{tabular}{@{}cc@{}}
\toprule
QAP Instances & $m$ \\
\midrule
had12 & 144 \\
had14 & 196 \\
had16 & 256 \\
had18 & 324 \\
had20 & 400 \\
 \bottomrule
\end{tabular}
\begin{tabular}{@{}cc@{}}
\toprule
QAP Instances & $m$ \\
\midrule
rou12 & 144 \\
rou15 & 225 \\
rou20 & 400 \\
tai40a & 1600 \\
tai40b & 1600 \\ \bottomrule
\end{tabular}
\begin{tabular}{@{}cc@{}}
\toprule
TSP Instances & $m$ \\
\midrule
bays29 & 784 \\
bayg29 & 784 \\
berlin52 & 2601 \\
brazil58 & 3249 \\
dantzig42 & 1681 \\
 \bottomrule
\end{tabular}
\begin{tabular}{@{}cc@{}}
\toprule
TSP Instances & $m$ \\
\midrule
fri26 & 625 \\
gr17 & 256 \\
gr21 & 400 \\
gr24 & 529 \\
st70 & 4761 \\ \bottomrule
\end{tabular}}
\end{table}

QUBO matrices (in upper triangular format) representing the cost and constraint functions of these QAP and TSP instances used are made available\footnote{\url{https://github.com/mayoayodelefujitsu/QUBOs}}

\subsection{Performance Measure}
\label{sec:arpd}
We compare the performance of the DA using different methods of generating penalty weights. The performance measure used is the Average Relative Percentage Deviation (ARPD) defined in Eq.~(\ref{eq:arpd})

\begin{equation}
\label{eq:arpd}
ARPD = \frac{1}{r}\sum_{i=1}^{r}\left ( \frac{\text{DA}(\alpha)_i -\text{Optimal} }{\text{Optimal}} \right ) \times 100\;.
\end{equation}

Note that $\text{DA}(\alpha)_i$ is the best energy returned by the DA for the $i^{th}$ run using penalty weight set to $\alpha$ and $r$ is the number of runs. We set $r=20$ and the optimal value are obtained from QAPLIB and TSPLIB.

\subsection{Parameter Settings: DA Algorithm}
\label{sec:pda}

The parameter settings used in the DA and DA3 are shown in Table \ref{tb:param}. For the DA, the temperature is set to decrease from `Initial Temperature' to `Final Temperature' by a fraction `Temperature Decay'. Note $\delta_i = \max\left ( \delta_f, \delta_{i-1} *\left ( 1-\rho\right ) \right )$ where $\delta_i$ denotes the temperature at iteration $i$. In the DA3, temperature and offset increment related parameters are automatically set. Therefore, no manual setting is required for the DA3, these parameters are thus shown as `NA' for DA3 in Table \ref{tb:param}. The stopping criteria used in the DA (CPU implementation) is `number of iterations' but `time (in seconds)' is used in DA3. This is because the two stopping criteria allowed in the DA3 are time and target energy (fitness). The DA is executed for $m^2$ iterations while the time limit for DA3 is set to the ceiling of 3\% of $m$ in seconds  ($m$ is presented in Table \ref{tb:inst}). Each experiment is executed independently 20 times. The parameters are chosen based on preliminary experiments. Note that $\beta = VLM$ (Eq.~(\ref{eq:wvermaa})).

\begin{table}[htb]
\begin{center}
\caption{Parameter Settings \label{tb:param}}
\resizebox{0.6\textwidth}{!}{  
\begin{tabular}{ccc}
\toprule
\textbf{Parameter} & \textbf{DA}  & \textbf{DA3} \\
\toprule
Initial Temperature $\delta_0$ & 0.1$\beta$, $\beta$, 10$\beta$ & NA\\
Final Temperature $\delta_f$ & 1 &NA\\
Stopping Criteria& $m^2$  & $\left \lceil 0.03m \right \rceil$ sec\\
offset\_increase\_rate &  $\delta_0 \div m^2$ & NA\\
Number of Runs $r$ & 20 & 20 \\
Temperature Decay $\rho$ & 0.001 & NA\\
\bottomrule
\end{tabular}}
\end{center}
\end{table}

\section{Results}
\label{sec:res}
In this section, different methods of generating penalty weights are compared using the parameter settings presented in Section \ref{sec:pda}. Results using the CPU implementation of the first generation DA algorithm are presented. Further results relating to the third generation DA are also presented.

\subsection{Penalty Weights for TSP and QAP instances}
Table \ref{tb:newtspqap} shows the penalty weights derived for different TSP and QAP instances using the methods of generating penalty weights defined in Section \ref{sec:pre}. The smallest and valid penalty weights are highlighted in bold. Validity of the penalty weights are assessed in Section \ref{sec:cpu}. It should be noted that the methods were applied to QUBO matrices in upper triangular format.

\begin{table}[h!tb]
\centering
\caption{Penalty weights for QAP and TSP instances derived using different methods \label{tb:newtspqap}}
\resizebox{0.9\textwidth}{!}{
\begin{tabular}{@{}ccccccc@{}}
\toprule
\multirow{2}{*}{\begin{tabular}[c]{@{}c@{}}Problem\\ Category\end{tabular}} & \multirow{2}{*}{\begin{tabular}[c]{@{}c@{}}Instance\\ Name\end{tabular}} & \multicolumn{5}{c}{Penalty Weights} \\ \cmidrule(l){3-7} 
 &  & UB & MQC & VLM & MOMC & MOC \\ \midrule
QAP & had12 & 249,240 & 126 & 5,460 & 2,730 & \textbf{488} \\
QAP & had14 & 573,484 & 162 & 8,968 & 4,484 & \textbf{533} \\
QAP & had16 & 1,014,488 & 162 & 12,580 & 6,290 & \textbf{545} \\
QAP & had18 & 1,832,940 & 200 & 16,102 & 8,051 & \textbf{1,513} \\
QAP & had20 & 2,950,640 & 220 & 20,928 & 10,464 & \textbf{1,335} \\
QAP & rou12 & 40,734,756 & 19,602 & 874,944 & 437,472 & \textbf{34,531} \\
QAP & rou15 & 98,340,328 & 19,602 & 1,498,176 & 749,088 & \textbf{79,715} \\
QAP & rou20 & 346,044,384 & 19,602 & 2,569,174 & 1,284,587 & \textbf{123,342} \\
QAP & tai40a & 5,904,547,332 & 19,602 & 10,418,804 & 5,209,402 & \textbf{176,904} \\
QAP & tai40b & 1,767,388,016,312 & 32,656,592 & 4,524,144,275 & 2,262,072,138 & \textbf{56,133,309} \\ \midrule
TSP & bayg29 & 3,381,534 & \textbf{386} & 6,279 & 3,140 & 2,404 \\
TSP & bays29 & 4,259,764 & \textbf{509} & 8,593 & 4,297 & 3,003 \\
TSP & berlin52 & 74,165,126 & \textbf{1,716} & 55,515 & 27,758 & 27,148 \\
TSP & brazil58 & 379,655,572 & \textbf{8,700} & 288,552 & 144,276 & 55,557 \\
TSP & dantzig42 & 4,814,472 & \textbf{192} & 5,029 & 2,515 & 1,915 \\
TSP & fri26 & 1,455,150 & \textbf{280} & 4,833 & 2,417 & 1,616 \\
TSP & gr17 & 1,005,188 & \textbf{745} & 7,981 & 3,991 & 3,074 \\
TSP & gr21 & 2,666,064 & \textbf{865} & 11,160 & 5,580 & 2,853 \\
TSP & gr24 & 1,609,942 & \textbf{389} & 5,185 & 2,593 & 1,888 \\
TSP & st70 & 16,647,424 & \textbf{129} & 5,055 & 2,528 & 2,079 \\ \bottomrule
\end{tabular}}
\end{table}

\subsection{CPU implementation of the DA Algorithm: comparing methods of generating penalty weights}
\label{sec:cpu}
In this section, results derived using CPU implementation of the first generation DA algorithm are presented. Table \ref{tb:feasrate} shows the number of feasible solutions returned by the DA within the stopping criteria when different methods of generating penalty weights are used. 

\begin{table}[t]
\begin{center}
\caption{Number of DA runs that returned feasible solutions out of 20 runs using different methods of generating weights for QAP instances (left) and TSP instances (right). The same number of feasible solutions was obtained for QAP and TSP instances with different values of initial temperature ($\delta_0 = 0.1\beta/\beta/10\beta$) apart from `rou12' when MOC is used, the respective values derived using each temperature is therefore presented.\label{tb:feasrate}}
\resizebox{0.8\textwidth}{!}{
\begin{tabular}{ccccccc}
\toprule
\multirow{2}{*}{\begin{tabular}[c]{@{}c@{}}Instance \\ Name\end{tabular}} &  \multicolumn{5}{c}{Number of feasible runs} \\ \cline{2-6}
  & UB & MQC & VLM & MOMC & MOC \\\toprule
had12 & 20 & 0 & 20 & 20 & 20 \\ \midrule
had14 & 20 & 0 & 20 & 20 & 20 \\ \midrule
had16 & 20 & 0 & 20 & 20 & 20 \\ \midrule
had18 & 20 & 0 & 20 & 20 & 20 \\ \midrule
had20 & 20 & 0 & 20 & 20 & 20 \\ \midrule
rou12 & 20 & 0 & 20 & 20 & 13/14/14 \\ \midrule
rou15 & 20 & 0 & 20 & 20 & 20 \\ \midrule
rou20 & 20 & 0 & 20 & 20 & 20 \\ \midrule
tai40a & 20 & 0 & 20 & 20 & 20 \\ \midrule
tai40b & 20 & 0 & 20 & 20 & 20 \\  \bottomrule
\end{tabular}
\begin{tabular}{ccccccc}
\toprule
\multirow{2}{*}{\begin{tabular}[c]{@{}c@{}}Instance \\ Name\end{tabular}} &  \multicolumn{5}{c}{Number of feasible runs} \\ \cline{2-6}
  & UB & MQC & VLM & MOMC & MOC \\\toprule
bayg29 & 20 & 20 & 20 & 20 & 20 \\ \midrule
bays29 & 20 & 20 & 20 & 20 & 20 \\ \midrule
berlin52& 20 & 20 & 20 & 20 & 20 \\ \midrule
brazil58 &20 & 20 & 20 & 20 & 20 \\ \midrule
dantzig42 &20 & 20 & 20 & 20 & 20 \\ \midrule
fri26 & 20 & 20 & 20 & 20 & 20 \\ \midrule
gr17 &20 & 20 & 20 & 20 & 20 \\ \midrule
gr21 &20 & 20 & 20 & 20 & 20 \\ \midrule
gr24 &20 & 20 & 20 & 20 & 20 \\ \midrule
st70 & 20 & 20 & 20 & 20 & 20 \\ \bottomrule
\end{tabular}
}

\end{center}
\end{table}

Tables \ref{tb:qaparpd1}, \ref{tb:qaparpd2} and \ref{tb:qaparpd3} respectively present the ARPD derived by the DA at initial temperature set to $0.1\beta, \beta$ and $10\beta$ using different methods of generating penalty weights. The ARPD for MQC is not presented for QAP instances in any of the tables because no feasible solution was found.

UB, MQC and MOC present their best ARPD averaged across TSP instances when the temperature is set to $0.1\beta$ while MOMC and VLM present their best ARPD on TSP instances when temperature is set to $\beta$. UB presents its best performance on QAP instances when temperature is set to $0.1\beta$, and MOC, MOMC and VLM present the best ARPD averaged across all QAP instances when the temperature is set to $10\beta$ .

\begin{table}[h!tb]
\begin{center}
\caption{ARPD from Optimal on QAP (left) and TSP (right) instances where initial temperature $\delta_0 = 0.1 \times \beta$ \label{tb:qaparpd1}}
\resizebox{\textwidth}{!}{
\begin{tabular}{cccccc}
\toprule
\multirow{2}{*}{\begin{tabular}[c]{@{}c@{}}Instance \\ Name\end{tabular}} & \multirow{2}{*}{Optimal} & \multicolumn{4}{c}{ARPD} \\\cline{3-6} 
 &  & UB & VLM & MOMC & MOC \\\toprule
had12 & 1,652 & 14.15 & 12.98 & 11.98 & \textbf{6.40} \\
had14 & 2,724 & 16.20 & 14.85 & 13.86 & \textbf{6.28} \\
had16 & 3,720 & 12.23 & 13.63 & 10.76 & \textbf{5.50} \\
had18 & 5,358 & 11.97 & 11.24 & 9.25 & \textbf{6.35} \\
had20 & 6,922 & 12.46 & 12.15 & 8.99 & \textbf{6.25} \\
rou12 & 235,528 & 29.12 & 29.15 & 27.98 & \textbf{\textit{10.37}} \\
rou15 & 354,210 & 30.75 & 33.34 & 28.21 & \textbf{16.28} \\
rou20 & 725,522 & 24.04 & 25.96 & 20.42 & \textbf{14.35} \\
tai40a & 3,139,370 & 20.44 & 20.73 & 16.08 & \textbf{13.00} \\
tai40b & 637,250,948 & 77.76 & 76.51 & 52.92 & \textbf{11.73} \\\bottomrule
 & Avg & 24.91 & 25.05 & 20.05 & \textbf{9.65}\\\bottomrule
\end{tabular}
\begin{tabular}{ccccccc}
\toprule
\multirow{2}{*}{\begin{tabular}[c]{@{}c@{}}Instance \\ Name\end{tabular}} & \multirow{2}{*}{Optimal} & \multicolumn{5}{c}{ARPD} \\ \cline{3-7}
  &  & UB & MQC & VLM & MOMC & MOC \\\toprule
bayg29 & 1,610 & 189.59 & \textbf{52.94} & 180.69 & 125.19 & 114.98 \\
bays29 & 2,020 & 190.45 & \textbf{55.52} & 188.23 & 130.47 & 116.61 \\
berlin52 & 7,542 & 295.70 & \textbf{100.40} & 289.46 & 212.58 & 209.45 \\
brazil58 & 25,395 & 389.04 & \textbf{138.94} & 390.04 & 276.79 & 260.21 \\
dantzig42 & 699 & 340.26 & \textbf{95.05} & 335.62 & 225.04 & 224.17 \\
fri26 & 937 & 177.06 & \textbf{57.32} & 177.34 & 120.84 & 107.51 \\
gr17 & 2,085 & 112.06 & \textbf{29.67} & 107.56 & 84.75 & 66.97 \\
gr21 & 2,707 & 170.91 & \textbf{44.82} & 166.01 & 123.01 & 93.32 \\
gr24 & 1,272 & 166.45 & \textbf{52.37} & 160.64 & 114.54 & 102.29 \\
st70 & 675 & 444.14 & \textbf{124.52} & 419.62 & 330.44 & 325.83 \\\bottomrule
 & Avg & 247.57 & 75.16 & 241.52 & 174.37 & 162.13\\\bottomrule
\end{tabular}
}
\end{center}
\end{table}

\begin{table}[h!tb]
\begin{center}
\caption{ARPD from Optimal on QAP (left) and TSP (right) instances where initial temperature $\delta_0 = \beta$}
\label{tb:qaparpd2}
\resizebox{\textwidth}{!}{
\begin{tabular}{cccccc}
\toprule
\multirow{2}{*}{\begin{tabular}[c]{@{}c@{}}Instance \\ Name\end{tabular}} & \multirow{2}{*}{Optimal} & \multicolumn{4}{c}{ARPD} \\ \cline{3-6}
 &  & UB & VLM & MOMC & MOC \\\toprule
had12 & 1,652 & 15.25 & 7.65 & 8.33 & \textbf{6.54} \\
had14 & 2,724 & 15.26 & 9.37 & 9.76 & \textbf{6.43} \\
had16 & 3,720 & 13.27 & 8.13 & 8.75 & \textbf{5.41} \\
had18 & 5,358 & 11.40 & 7.08 & 7.04 & \textbf{6.55} \\
had20 & 6,922 & 12.86 & 7.38 & 7.66 & \textbf{6.74} \\
rou12 & 235,528 & 32.12 & 20.34 & 20.75 & \textbf{\textit{9.58}} \\
rou15 & 354,210 & 31.33 & 22.00 & 21.37 & \textbf{15.75} \\
rou20 & 725,522 & 24.49 & 17.77 & 17.91 & \textbf{13.69} \\
tai40a & 3,139,370 & 20.43 & 16.10 & 16.13 & \textbf{12.89} \\
tai40b & 637,250,948& 78.61 & 51.81 & 51.25 & \textbf{11.49} \\\bottomrule
 & Avg & 25.50 & 16.76 & 16.89 & \textbf{9.51} \\\bottomrule
\end{tabular}

\begin{tabular}{ccccccc}
\toprule
\multirow{2}{*}{\begin{tabular}[c]{@{}c@{}}Instance \\ Name\end{tabular}} & \multirow{2}{*}{Optimal} & \multicolumn{5}{c}{ARPD} \\ \cline{3-7}
  &  & UB & MQC & VLM & MOMC & MOC \\\toprule
 bayg29 & 1,610 & 194.42 & \textbf{54.69} & 127.05 & 120.22 & 117.59 \\
bays29 & 2,020 & 196.69 & \textbf{57.22} & 124.98 & 122.92 & 120.10\\
berlin52 & 7,542 & 298.15 & \textbf{102.47} & 217.57 & 214.73 & 214.16 \\
brazil58 & 25,395 & 380.88 & \textbf{137.75} & 278.99 & 277.66 & 261.90 \\
dantzig42 & 699 & 350.45 & \textbf{100.25} & 238.07 & 232.41 & 222.18 \\
fri26 & 937 & 185.14 & \textbf{59.04} & 116.29 & 114.82 & 109.82 \\
gr17 & 2,085 & 128.08 & \textbf{31.41} & 70.44 & 62.56 & 60.52 \\
gr21 & 2,707 & 190.91 & \textbf{52.99} & 114.31 & 105.63 & 98.48 \\
gr24 & 1,272 & 178.25 & \textbf{52.85} & 114.77 & 104.18 & 98.75 \\
st70 & 675 & 435.61 & \textbf{129.66} & 335.66 & 331.78 & 329.34 \\\bottomrule
 & Avg & 253.86 & \textbf{77.83} & 173.81 & 168.69 & 163.28 \\\bottomrule
\end{tabular}}
\end{center}
\end{table}

\begin{table}[h!tb]
\begin{center}
\caption{ARPD from Optimal on QAP (left) and TSP (right) instances where initial temperature $\delta_0 = 10 \times \beta$\label{tb:qaparpd3}}
\resizebox{\textwidth}{!}{
\begin{tabular}{cccccc}
\toprule
\multirow{2}{*}{\begin{tabular}[c]{@{}c@{}}Instance \\ Name\end{tabular}} & \multirow{2}{*}{Optimal} & \multicolumn{4}{c}{ARPD} \\\cline{3-6}
 &  & UB & VLM & MOMC & MOC \\\toprule
had12 & 1,652 & 11.26 & 7.99 & 8.51 & \textbf{6.22} \\
had14 & 2,724 & 15.56 & 9.13 & 9.48 & \textbf{6.11} \\
had16 & 3,720 & 14.02 & 8.19 & 8.19 & \textbf{5.12} \\
had18 & 5,358 & 11.80 & 7.07 & 7.31 & \textbf{6.03} \\
had20 & 6,922 & 12.57 & 7.32 & 7.33 & \textbf{6.43} \\
rou12 & 235,528 & 28.30 & 18.94 & 16.50 & \textbf{\textit{10.02}} \\
rou15 & 354,210 & 33.98 & 21.02 & 20.16 & \textbf{14.57} \\
rou20 & 725,522 & 25.19 & 17.80 & 17.36 & \textbf{13.05} \\
tai40a & 3,139,370 & 20.96 & 16.00 & 15.97 & \textbf{12.54} \\
tai40b & 637,250,948 & 79.65 & 50.85 & 49.94 & \textbf{12.10} \\\bottomrule
 & Avg & 25.33 & 16.43 & 16.07 & \textbf{9.22} \\\bottomrule
\end{tabular}

\begin{tabular}{ccccccc}
\toprule
\multirow{2}{*}{\begin{tabular}[c]{@{}c@{}}Instance \\ Name\end{tabular}} & \multirow{2}{*}{Optimal} & \multicolumn{5}{c}{ARPD} \\ \cline{3-7}
  &  & UB & MQC & VLM & MOMC & MOC \\\toprule
bayg29 & 1,610 & 193.29 & \textbf{57.27} & 127.64 & 122.24 & 122.83 \\
bays29 & 2,020 & 196.84 & \textbf{57.57} & 132.84 & 124.92 & 111.88 \\
berlin52 & 7,542 & 300.65 & \textbf{103.48} & 218.34 & 215.59 & 214.92 \\
brazil58 & 25,395 & 375.04 & \textbf{139.30} & 285.40 & 273.91 & 265.80 \\
dantzig42 & 699 & 333.03 & \textbf{100.39} & 238.24 & 227.12 & 230.27 \\
fri26 & 937 & 180.85 & \textbf{62.51} & 124.09 & 114.74 & 11.27 \\
gr17 & 2,085 & 122.95 & \textbf{30.19} & 70.65 & 64.84 & 60.17 \\
gr21 & 2,707 & 178.13 & \textbf{52.73} & 115.05 & 107.30 & 99.54 \\
gr24 & 1,272 & 179.79 & \textbf{56.82} & 116.19 & 106.69 & 103.71 \\
st70 & 675 & 452.83 & \textbf{126.34} & 337.61 & 334.76 & 325.69 \\ \bottomrule
 & Avg & 251.34 & \textbf{78.66} & 176.61 & 169.21 & 164.56 \\\bottomrule
\end{tabular}}
\end{center}
\end{table}

In Table \ref{tb:qaparpd1}, the DA presents the best average ARPD on QAP instances when the method of setting penalty weight is set to MOC. The DA however presents the best average ARPD on TSP instances when MQC is used. These results are expected since the MOC and MQC respective present the smallest yet valid penalty weights for QAP and TSP instances. ARPD for QAP instance \textit{rou12}, is shown in italics when the method of setting penalty weight is set to MOC because the ARPD is only computed using the 13 feasible solutions found while others are generated using 20 feasible solutions. The DA presents the worst average ARPD on QAP (or TSP) instances when the method of setting penalty weights is set to VLM  (or UB).

In Tables \ref{tb:qaparpd2} and \ref{tb:qaparpd3}, similar to the results produced in Table \ref{tb:qaparpd1}, the DA presents the best average ARPD on QAP instances when the methods of setting penalty weight is set to MOC and the best average ARPD on TSP instances when set to MQC. ARPD for QAP instance \textit{rou12}, is shown in italics when the method of setting penalty weight is set to MOC because the ARPD is only computed using the 14 feasible solutions found while others are generated using 20 feasible solutions. The DA presents the worst average ARPD on QAP and TSP instances when the method of setting penalty weights is set to UB.

In general, the results show that the methods which produced the smallest valid penalty weights (MQC for TSP and MOC for QAP) consistently produced the best ARPD. Conversely, the results also show that the method that produced the largest penalty weights on TSP and QAP instances (UB) often presents the worst ARPD. For permutation problems represented as two-way one-hot, all neighbours of a feasible solution are infeasible solutions. It is therefore important for penalty weights to be small enough to encourage the algorithm to explore infeasible solutions in order to find better feasible solutions.

\subsection{DA3: Comparing methods of generating penalty weights}
In section \ref{sec:cpu}, we show how different methods of generating penalty weights can affect the performance of the 1st generation DA (CPU implementation). In this section, we present results using DA3 (i.e. version of the third generation DA which benefits from hardware speed-up). It should be noted that the DA3 has more capabilities and can be executed in many modes. A major improvement presented by the DA3 is the ability to handle linear inequality and one-hot constraints.

\begin{table}[h]
\centering
\caption{ARPD derived using DA3 with different methods of generating penalty weights on QAP (left) and TSP (right) instances\label{tb:arpdda3}}
\resizebox{0.92\textwidth}{!}{
\begin{tabular}{@{}cccccc@{}}
\toprule
\multirow{2}{*}{\begin{tabular}[c]{@{}c@{}}Instance \\ Name\end{tabular}} & \multirow{2}{*}{Optimal} & \multicolumn{4}{c}{ARPD} \\ \cmidrule(l){3-6} 
 &  & UB & VLM & MOMC & MOC \\ \midrule
had12 & 1,652 & \textbf{0.00} & \textbf{0.00} & \textbf{0.00} & \textbf{0.00} \\
had14 & 2,724 & \textbf{0.00} & \textbf{0.00} & \textbf{0.00} & \textbf{0.00} \\
had16 & 3,720 & \textbf{0.00} & \textbf{0.00} & \textbf{0.00} & \textbf{0.00} \\
had18 & 5,358 & \textbf{0.00} & \textbf{0.00} & \textbf{0.00} & \textbf{0.00} \\
had20 & 6,922 & \textbf{0.00} & \textbf{0.00} & \textbf{0.00} & \textbf{0.00} \\
rou12 & 235,528 & \textbf{0.00} & \textbf{0.00} & \textbf{0.00} & \textbf{0.00} \\
rou15 & 354,210 & \textbf{0.00} & \textbf{0.00} & \textbf{0.00} & \textbf{0.00} \\
rou20 & 725,522 & \textbf{0.00} & \textbf{0.00} & \textbf{0.00} & \textbf{0.00} \\
tai40a & 3,139,370 & \textbf{0.07} & \textbf{0.07} & \textbf{0.07} & \textbf{0.07} \\
tai40b & 637,250,948 & \textbf{0.00} & \textbf{0.00} & \textbf{0.00} & \textbf{0.00} \\ \midrule
\multicolumn{2}{c}{Average} & 0.01 & 0.01 & 0.01 & 0.01 \\ \bottomrule
\end{tabular}
\begin{tabular}{@{}ccccccc@{}}
\toprule
\multirow{2}{*}{\begin{tabular}[c]{@{}c@{}}Instance \\ Name\end{tabular}}  & \multirow{2}{*}{Optimal} & \multicolumn{5}{c}{ARPD} \\ \cmidrule(l){3-7} 
 &  & UB & MQC & VLM & MOMC & MOC \\ \midrule
bayg29 & 1,610 & \textbf{0.00} & \textbf{0.00} & \textbf{0.00} & \textbf{0.00} & \textbf{0.00} \\
bays29 & 2,020 & \textbf{0.00} & \textbf{0.00} & \textbf{0.00} & \textbf{0.00} & \textbf{0.00} \\
berlin52 & 7,542 & \textbf{1.86} & 2.96 & 4.16 & 6.05 & 2.20 \\
brazil58 & 25,395 & 1.53 & 3.33 & 1.53 & 1.58 & \textbf{1.43} \\
dantzig42 & 699 & \textbf{0.00} & \textbf{0.00} & \textbf{0.00} & \textbf{0.00} & \textbf{0.00} \\
fri26 & 937 & \textbf{0.00} & \textbf{0.00} & \textbf{0.00} & \textbf{0.00} & \textbf{0.00} \\
gr17 & 2,085 & \textbf{0.00} & \textbf{0.00} & \textbf{0.00} & \textbf{0.00} & \textbf{0.00} \\
gr21 & 2,707 & \textbf{0.00} & \textbf{0.00} & \textbf{0.00} & \textbf{0.00} & \textbf{0.00} \\
gr24 & 1,272 & \textbf{0.00} & \textbf{0.00} & \textbf{0.00} & \textbf{0.00} & \textbf{0.00} \\
st70 & 675 & 2.67 & \textbf{1.48} & 2.41 & 2.37 & 2.39 \\ \midrule
\multicolumn{2}{c}{Average} & 0.61 & 0.78 & 0.81 & 1.00 & 0.60 \\ \bottomrule
\end{tabular}}

\end{table}

We present the ARPD achieved within $0.03m$ seconds (where $m$ represents the size of the solution) of executing DA3 with different methods of generating penalty weights in Table \ref{tb:arpdda3}. DA3 achieves 100\% feasibility on TSP instances with any of the penalty methods, it also achieves 100\% feasibility on QAP when UB, VLM, MOMC and MOC methods are used. Furthermore, the standard deviation across 20 runs of the DA3 is often 0. For all QAP instances as well as 7 out of 10 TSP instances, DA3 obtains the same ARPD regardless of the penalty weights. DA3 presents varying ARPD on the largest TSP instances when different methods of generating penalty weights are used. Best ARPD was obtained for berlin52, brazil58 or st70 when UB, MOC or MQC is used respectively. There is therefore no clear evidence of one method being the best. We can therefore not make the same conclusions as the previous section, where there was clear evidence of smaller and valid penalty weights leading to better solution quality. Similarities in performance of DA3 regardless of penalty weights used may be because of the algorithmic changes made since 1st generation DA. DA3 is designed to solve problems formulated as two-way one-hot more efficiently \cite{da3}. It is also able to automatically find the best parameter settings for any BQP. The algorithm that supports DA3 is however not publicly available, it is therefore difficult to be precise about the reason for the difference in performance when compared to the first generation DA.

\section{Conclusion and Further Work}
\label{sec:con}
Permutation problems like TSP and QAP can be formulated as QUBO such that algorithms that use specialised hardware e.g. Quantum Annealer or DA can solve them. Transforming these problems into QUBO form requires the setting of penalty weights. In this study, we examined different methods of generating penalty weights within the context of using the DA algorithm for solving permutation problems. The permutation problems used are TSP and QAP. We present improvements to existing methods of generating penalty weights leading to better performance of the DA. Although the DA algorithm, which shares similar properties with SA, was influenced by the magnitude of penalty weights, we could not reach the same conclusions with DA3. It was impossible to do deeper analysis because the DA3 algorithm is not publicly available. Further research into how various algorithms behave with different mechanisms of generating penalty weights is therefore necessary.